\input amstex
\documentstyle{amsppt}
\NoBlackBoxes
\magnification=1200
\parindent 20 pt
\vsize=7.50in

\define\1{^{-1}}

\define\RP{\Bbb R \Bbb P^2}

\define \C{\Bbb C}
\define \Z{\Bbb Z}
\define \Q{\Bbb Q}
\define \R{\Bbb R}
\define \p{\Bbb P}
\define \a{\alpha }
\define \li{L}

\define \bk{\bigskip}
\define \ov{}
\define \edm{\enddemo}
\define \ep{\endproclaim}
\define\pr{\operatorname{pr}}

{\catcode`\@=11
\gdef\n@te#1#2{\leavevmode\vadjust{%
 {\setbox\z@\hbox to\z@{\strut#1}%
  \setbox\z@\hbox{\raise\dp\strutbox\box\z@}\ht\z@=\z@\dp\z@=\z@%
  #2\box\z@}}}
\gdef\leftnote#1{\n@te{\hss#1\quad}{}}
\gdef\rightnote#1{\n@te{\quad\kern-\leftskip#1\hss}{\moveright\hsize}}
\gdef\?{\FN@\qumark}
\gdef\qumark{\ifx\next"\DN@"##1"{\leftnote{\rm##1}}\else
 \DN@{\leftnote{\rm??}}\fi{\rm??}\next@}}

\def\Aut{\operatorname{Aut}}
\def\Kl{\operatorname{Kl}}
\def\kl{\operatorname{kl}}
\def\div{\operatorname{div}}
\def\tr{\operatorname{tr}}

\def\Div{\operatorname{Div}}
\def\Pic{\operatorname{Pic}}
\def\Id{\operatorname{Id}}

\def\+{\sqcup}
\def\*{!}

\topmatter
\title
On real structures of rigid surfaces
\endtitle
\author V. Kharlamov and Vik. S. Kulikov$^*$
\endauthor
\abstract
We construct several rigid (i.e., unique in their
deformation class) surfaces which
have particular behavior
with respect to real structures: in one example
the surface has no any real structure, in the other
one it has a unique real structure and this structure is not maximal
with respect to the Smith-Thom inequality.
So, it answers in negative
to the following problems:
existence of real surfaces
in each complex deformation class and existence of
maximal surfaces in each complex deformation class containing
real surfaces.
Besides, we prove that there is no real surfaces among the surfaces of general
type  with $p_g=q=0$ and $K^2=9$.

The surfaces constructed provide new counter-examples to the
``Dif=Def'' problem.
\endabstract
\thanks \flushpar $^*$ This work was done during the stay of the second
author in Strasbourg university and it is partially supported by
INTAS-OPEN-97-2072, NWO-RFBR-047-008-005, and
RFBR  99-01-01133.\endthanks
\endtopmatter

\document

\baselineskip 20pt

\subheading{\S 0.\ Introduction}

One of the principal settings in real algebraic geometry
is to fix a deformation class of complex
varieties and to study, inside this class,
the varieties which can be equipped with a real
structure (and then investigate their topological,
as well as
other invariant under real deformations, properties).
Those, which are maximal with respect to
the Smith-Thom bound, are of a special interest
(since they have spectacular topological properties,
see, for example, the survey \cite{DK}; for surfaces
this bound is reproduced below in Section 5).
Thus, two natural questions arise: does any complex
deformation class of compact complex varieties contain
a real variety; and does any complex deformation class
containing real varieties contain a maximal one?
Up to our knowledge,
in dimension $\ge 2$
the both questions remained
open till now.
We show that the response
to the both questions is in negative.
In our examples
the varieties are surfaces which are rigid,
where the latter means that the moduli space
of complex structures on the underlying smooth manifold is
$0$-dimensional.
Moreover, in our examples they are strongly rigid,
i.e.,
the quotient of the moduli space
by the canonical complex conjugation
(which replaces a complex structure
of the surface
by the complex conjugated
one, and thus holomophic functions by anti-holomophic ones;
the orientation of the underlying smooth 4-manifold is preserved)
is merely a point. It is worth noticing that
in the first of our
examples the moduli space consists of two conjugated points, in the
second one it reduces to one real point, see the remarks in Section 4.
Besides,
the
two conjugated surfaces in the first example give
one more counterexample to "Dif=Def" problem
(earliest counterexamples
were constructed by Manetti in \cite {Ma}).
In fact, in
all our examples
the surfaces are of general type
and with $c_1^2=3c_2$ (they are the so-called
Miyaoka-Yau surfaces; the fact that they are strongly rigid
and, moreover, unique in their homotopy
type up to
holomorphic and anti-holomorphic diffeomorphisms
is well known, see, for example,
\cite{BPV}).
Following F.~Hirzebruch
\cite{H}
we construct such rigid
surfaces as (finite abelian) Galois
coverings of the (blown-up)
projective plane
branched along arrangements of lines.
We start from giving in Sections 1 and 2
their explicit construction
via the orbit spaces of the Ferma covering.
In Section 3 we study the group of automorphisms and
anti-automorphisms of the constructed surfaces. In Section 4
three
main examples are treated.
In Section 5
we consider fake projective planes
(that is, the surfaces of general type with $c_2=3$ and $c_1^2=9$).
We prove that they have no
anti-holomorphic
diffeomorphisms and,
in particular,
can not be equiped
with a real structure.
This section contains also several
remarks on
other related topics.

{\bf Acknowledgements.} We are grateful to Y.~Miyaoka for
stimulating our interest
to the real geometry of rigid surfaces
and T.~Delzant
for useful proposals during
the preparation of this publication. The
first author is grateful to F.~Catanese
for an intersting discussion
of the deformation
problems in real algebraic geometry.
\bk

\subheading{\S 1. Galois coverings of the plane
branched over an arrangement of lines\ }

By a Galois covering
of a smooth algebraic variety $Y$
we mean a finite morphism $f:X\to Y$ of
a normal algebraic variety $X$ to $Y$
such that the function fields
imbedding $\C (Y)\subset \C (X)$
induced by $f$ is a Galois extension.
As is
well known, a finite morphism
$f:X\to Y$ is a Galois covering with
Galois group $G$ if and only if
$G$ coincides with the group of covering transformations
and the latter acts transitively
on every fiber of $f$.
Besides, a finite branched covering
is Galois if and only if the unramified
part of the covering
(i.e., the restriction to the complements of the ramification and branch
loci)
is Galois. In addition, a branched covering
is determined up to isomorphism by its unramified part
and, moreover,
a covering morphism from the unramified part
of one branched covering to the unramified part of another one
induces a covering morphism between these branched coverings
if the extension of the morphism of underlying varieties to the branch
loci is given.
Let
us recall also that an unramified covering is Galois with Galois group $G$
if and only if
it is a covering associated with an epimorphism of the fundamental group
of the underlying variety to $G$, and, in particular,
the Galois coverings with abelian Galois group $G$ are in one-to-one
correspondence with epimorphisms to $G$ of the first homology group with
integral coefficients.
All these results are well known and their most
nontrivial part can be deduced,
for example, on the Grauert-Remmert existence theorem
\cite{G-R} (a detailed exposition of the basic results on
branched coverings is found, f.e., in \cite{N}).

In what follows we have deal only with
coverings of the complex
projective plane $\p ^2$
ramified over an arrangement of lines
$\li=\li _1\cup \dots \cup \li _n$.
Similarly to general abelian Galois coverings,
a Galois covering $g:Y \to \p ^2$ of
$\p ^2$ with abelian Galois group $G$ branched along
$\li $ is determined
uniquely by an epimorphism $\ov{\varphi} : H_1(\p ^2 \setminus \li ,
\Z )\to G$, and it exists for any such an epimorphism.
Since $H_1(\p ^2 \setminus \li ,\Z )\simeq \Z ^{n-1}$,
there exists, in particular, a covering
$g_{u(m)}: Y_{u(m)}\to \p ^2$ corresponding to the natural epimorphism
$\ov{\varphi }:  H_1(\p ^2 \setminus \li , \Z )\to H_1(\p ^2
\setminus \li ,\Z/m\Z )= H_1(\p ^2 \setminus \li ,\Z )\otimes (\Z
/m\Z )$\,. We call it {\it Ferma covering}\,.
The following statement is also an immediate consequence of
the general results on branched coverings recalled
in the beginning of this Section.

\proclaim{Proposition 1.0}
If $g: Y \to \p ^2$ is a Galois covering with Galois
group $G\simeq (\Z /m\Z )^k$ branched along $\li $, then
$k\le n-1$ and for any epimorphism
$H_1(\p^2\setminus\li)\to G$
there exists a unique Galois covering
$f: Y_{u(m)}\to Y$
inducing this epimorphism and such that
$g_{u(m)}=g\circ f$\,.
\qed
\endproclaim

In what follows we have deal with
Galois coverings whose Galois group is $G\simeq (\Z /m\Z )^k$,
and we construct them in a way described in
the above proposition.

The simple loops
$\lambda _i, 1\le i\le n,$ around the lines $\li _i$
generate $H_1(\p ^2 \setminus \li ,\Z )\simeq \Z ^{n-1}$.
They are subject to the relation
$$\lambda _1+\dots +\lambda _n=0,$$
and without
loss of generality we can assume that the universal covering
$g_{u(m)}: Y_{u(m)}\to \p ^2$ is determined by
the epimorphism $\ov{\varphi }: H_1(\p ^2 \setminus \li ,
\Z )\to (\Z /m\Z )^{n-1}$ sending
$\lambda _n$ to $(m-1,\dots , m-1)$ and
$\lambda _i$ with $1\leq i\leq n-1$
to $(0,\dots,0,1,0,\dots,0)$
with $1$ in the $i$-th place.
We choose an additional line
$\li _{\infty}\subset \p ^2$ in general
position with respect to $\li $ and
introduce affine coordinates $(x_1,x_2)$ in
$\C ^2=\p ^2\setminus \li _{\infty}$. Let $l_i(x_1,x_2)=0$ be
a linear
equation of $\li _i \cap \C ^2$. Put $z_i=(l_il_{n}^{m-1})^{1/m}$,
$1\leq i \leq n-1$.
Then the function field $K_{u(m)}=\C
(Y_{u(m)})=\C (x_1,x_2,z_1,\dots ,z_{n-1})$ of
$Y_{u(m)}$ is the abelian extension of the function field
$k=\C (x_1,x_2)$ of $\p ^2$ of degree $m^{n-1}$
with Galois group
$$G=\{ \, \ov{\gamma}=
(\gamma_1,\dots ,\gamma _n)\in (\Z /m\Z )^{n}\, \mid \, \sum \gamma _i \equiv
0 \, (\text{mod}\, m)\, \}\simeq (\Z /m\Z )^{n-1}.$$
(In other words, the pull-back of
$\p^2\setminus \li_\infty$ in $Y_{u(m)}$ is
naturally isomorphic to the normalization of the
affine
subvariety of $\C^{n+1}$ given in
coordinates $x_1,x_2,z_1,\dots,
z_{n-1}$ by equations $z_1^m=l_1l_n^{m-1},\dots z_{n-1}^m=l_{n-1}l_n^{m-1}$.)

For a multi-index
$\ov{a}=(\alpha _1,\dots ,\alpha_{n-1})$,
 $0\leq \alpha _i\leq m-1$, we put
$$z^{a}=\prod_{i=1}^{n-1}z_i^{\a _i}.$$
The action of $\ov{\gamma}=
(\gamma_1,\dots ,\gamma _n)\in G$ on $K_{u(m)}$ is given by
$$\ov{\gamma}(z^{\ov{a}})=\mu ^{(\ov{\gamma},\ov{a})}
z^{\ov{a}},$$
where
$$(\ov{\gamma},\ov{a})=\sum_{j=1}^{n-1}\gamma _j\a _j$$
and $\mu =e^{2\pi i/m}$ is the $m$-th root of the unity.
Thus,
$$K_{u(m)}=\bigoplus_{
0\le \a _i\le m-1}\C (x_1,x_2)z^{\ov{a}}$$
is a decomposition of the vector space $K_{u(m)}$ over $\C (x_1,x_2)$
into a finite direct sum of
degree $1$ representations of $G$.

Let $\ov{\varphi} : H_1(\p ^2 \setminus \li ,\Z )\to (\Z /m\Z)^k$ be
an epimorphism given by $\ov{\varphi}(\lambda _i)=(\a _{i,1},\dots ,\a _{i,k})$,
where $\a _{1,j}+\dots +\a _{n,j}\equiv 0\, \text{mod}\, m$ for
every $j=1,\dots , k$, and let
$g:Y \to \p ^2$ be the corresponding Galois covering.
Then, by Proposition 1.0, there exists
a unique Galois covering $f: Y_{u(m)}\to Y$.
It determines the inclusion
$f^{*} :\C (Y)\to K_{u(m)}$
of the function field $\C (Y)$ of $Y$
into the function field $K_{u(m)}=\C(Y_{u(m)})$.
Clearly, $\C(Y)$ is the subfield
$K_{\ov{\varphi}}=\C (x_1,x_2,w_1,\dots ,w_{k})$ of $K_{u(m)}$,
where
$w_j=z_1^{\a _{1,j}}\cdot \dots \cdot 
z_{n-1}^{a_{n-1,j}}$, and
$$\text{Gal}(K_{u(m)}/K_{\ov{\varphi}})= \{ \,
(\gamma_1,\dots ,\gamma _n)\in G\, \mid \, \sum_{i=1}^{n-1} \a _{i,j}\gamma_{i}\equiv 0\,
\text{mod}\, m,\, 1\leq j \leq k\,  \}.$$

By construction,
$Y$ is a normal surface with isolated singularities.
The singular points of $Y$ can appear only
over the $r$-fold points of $\li$ with $r\ge 2$,
i.e., over points lying on $r$ lines $\li _{i_1},\dots ,\li _{i_r}$
of the arrangement.

In what follows
we call $r$ elements of $(\Z/m\Z)^k$ {\it linear independent
over} $\Z/m\Z$ if they generate in $(\Z/m\Z)^k$
a subgroup isomorphic to $(\Z/m\Z)^r$ (and thus
admitting $(\Z/m\Z)^{k-r}$
as
its complement).

\proclaim{Lemma 1.1}
Let $p$ be a $2$-fold point of $\li$ and $\ov{\varphi}(\lambda_{i_1})$ and
$\ov{\varphi}(\lambda_{i_2})$ are linear independent over $\Z /m\Z $
in $(\Z /m\Z )^k$.
Then the surface $Y$ is non-singular at
each point of $f^{-1}(p)$.
\ep

{\it Proof}. Let $p= L_{i_1}\cap L_{i_2}$.
Choose a small round neighborhood $U$
of $p$ in $\p ^2$ and local analytic coordinates
$y_1,y_2$  in $U$
such that $y_j=0$ is an equation of $L_{i_j}$. Then,
$H_1(U\setminus (L_{i_1}\cup L_{i_2}),\Z )\simeq \Z\oplus \Z$.
At any point $q\in g^{-1}(p)$ the germ $V\to U$
of the covering $Y\to\p^2$ is a $G'$-covering, where
$G'$ is the image of $H_1(U\setminus (L_{i_1}\cup L_{i_2}),\Z)$
under the composition $\varphi \circ i_*$ of $\varphi$
with the inclusion homomorphism $i_* : H_1(U\setminus (L_{i_1}\cup L_{i_2}),\Z )\to
H_1(\p ^2\setminus L,\Z )$. Moreover, this $G'$-covering is determined
by $\varphi\circ i_*$. Identifying
$\ov{\varphi}(\lambda _{i_1})$,
$\ov{\varphi}(\lambda _{i_2})$
with the standard generators of $(\Z/m\Z)^2$
we get an isomorphism between $V\to U$
and the covering determined by equations
$z_1^m=y_1, z_2^m=y_2$. Thus, $V$ is nonsingular.
\qed

In our further examples,
to resolve the singularities of $Y$
over the $r$-fold points of $L$
with $r\geq 3$, we blow up all these points.
Let $\sigma :\widetilde {\p ^2} \to \p ^2$
be this blow up,
$L'_i$ the strict transform of $\li _i$,
$E_p$ the line blown up over a $r$-fold point $p$, and
$\varepsilon_{p}\in H_1(\widetilde{\p ^2}\setminus \sigma ^{-1}(\li ),\Z )
=H_1(\p ^2 \setminus \li ,\Z )$
a simple loop around $E_{p}$.

The identification $H_1(\widetilde{\p ^2}
\setminus \sigma ^{-1}(\li ),\Z )
=H_1(\p ^2 \setminus \li ,\Z )$
composed with $\varphi$ provides
an epimorphism
$\varphi : H_1(\widetilde{\p ^2}
\setminus \sigma ^{-1}(\li ),\Z )\to (\Z /m\Z )^k$.
Let consider the associated Galois covering
$f:X \to \widetilde{\p ^2}$.

The proof of the following statements is straightforward
(to establish the relation given by the first statement it is sufficient
to consider a generic line pencil; the second statement follows
from Lemma 1.1).

\proclaim{Lemma
1.2} Let $p=\li _{i_1}\cap \dots \cap \li _{i_r}$ be
an $r$-fold point of $\li$. Then $\varepsilon_{p}=\lambda _{i_1}+\dots +
\lambda _{i_r}$.\qed
\ep

\proclaim{Lemma
1.3} If for each $r$-fold point
$p=\li _{i_1}\cap \dots \cap \li _{i_r}$ of $\li$ with $r\geq 3$
the pairs $\varphi(\varepsilon_{p})$ and $\varphi(\lambda _{i_j})$,
$1\leq j \leq r$, are linear independent over $\Z /m\Z $ in $(\Z /m\Z )^k$,
then $X$ is nonsingular.
\qed
\ep

As a consequence,
the constructed surface $X$ is a resolution of singularities
of $Y$. Indeed, the covering $f$
is included in the commutative diagram
$$
\CD
X @>{
\pi }>>
Y \\
@V{f}VV            @VV{
g}V     \\
\widetilde{\p ^2} @>>{\sigma}>  \p ^2.
\endCD
$$
with a regular map $\pi$
(clearly, it is continuous, and thus
its regularity follows,
for example, from regularity
on $X\setminus f^{-1}(\sigma^{-1}(\li))$).

\subheading{\S 2. $(\Z /5\Z )^2$-Galois coverings
branched over the arrangement of lines dual
to the inflection points of a
smooth
cubic\ }

We use the notation of \S 1.

Let $\li=\li _1\cup \dots \cup \li _9$ be  an arrangement of nine
lines in $\p ^2$ dual to the nine inflection points of a smooth cubic $C$
in the dual plane.
Let $t_r$ ($r\geq 2$) be the number of $r$-fold points of $\li$,
i.e., the number of points lying on exactly $r$ lines of the arrangement.
As is well-known (and easy to check using the group law on the cubic),
in this arrangement
$t_3=12$, $t_r=0$ if $r\neq 3$,
and exactly four singular points of $\li$ lie on each $\li_i$,
$1\le i\le 9$.
(Note that the arrangement of lines dual to the inflection
points of a smooth cubic is rigid, i.e., any such arrangement can be
transformed to another by a linear transformation of the projective plane.)

If $C$ is a cubic given by
$x_1^3+x_2^3+x_3^3=0$,
then the lines
$\li _1, \dots \li _9$ are given by equations
$$
\matrix
\li _1=\{x_1-x_3=0 \}, & \li _2=\{x_1-\mu^2 x_3=0 \}, &
\li _3=\{x_1+\mu x_3=0 \}, \\
\li _4=\{x_2-\mu^2x_3=0 \}, & \li _5=\{x_2-x_3=0 \}, &
\li _6=\{x_2+\mu x_3=0  \},  \\
\li _7=\{x_1+\mu x_2=0  \}, & \li _8=\{x_1-\mu^2x_2=0 \}, &
\li _9=\{ x_1-x_2=0 \}.
\endmatrix
$$
where $\mu =e^{\pi i/3}$.

The intersection of three distinct lines $\li_i,\li_j,\li_k$
is nonempty if and only if $(i,j,k)\in T$, where
$$
\align
T= & \{ (1,2,3), (4,5,6), (7,8,9), (1,4,7), (2,5,8), (3,6,9), \\
 & (1,5,9), (3,5,7), (1,6,8),
 (3,4,8), (2,4,9), (2,6,7) \}.
\endalign
$$
Denote by $p_{i,j,k}$, $(i,j,k)\in T$,
the point of intersection of $\li_i,\li_j,\li_k$.

Consider a Galois covering
$g:Y \to \p ^2$ with Galois group $G\simeq
(\Z /5\Z)^2$ branched along
$\li $ and  determined by an epimorphism $\ov{\varphi }:
H_1(\p ^2 \setminus \li , \Z )\to G$.

Denote by $\sigma :\widetilde {\p ^2} \to \p ^2$ the blow up
with the centers at all the $3$-fold points $p_{i,j,k}, (i,j,k)\in T$,
by $E_{i,j,k}$ the exceptional divisor over $p_{i,j,k}$,
and by $L'_i$ the strict transform
of $\li _i$. Let $\varepsilon_{i,j,k}\in
H_1(\widetilde{\p ^2}\setminus \sigma ^{-1}(\li ),\Z )\simeq
H_1(\p ^2 \setminus \li ,\Z )$
correspond to a simple loop around $E_{i,j,k}$.

The epimorphism $\varphi : H_1(\widetilde{\p} ^2 \setminus \sigma^{-1}
(\li ),\Z )\to (\Z /5\Z)^2$ determines a Galois covering $f:X\to
\widetilde{\p} ^2$.
Pose $C_i=f^{-1}(L'_i)$ and $D_{i,j,k}=f^{-1}(E_{i,j,k})$.

In what follows we assume that the epimorphism
$\varphi : H_1(\widetilde{\p} ^2 \setminus \sigma^{-1}
(\li ),\Z )\to (\Z /5\Z)^2$ satisfies the following condition

{\bf (S)} {\it $\varphi(\varepsilon_{i_1,i_2,i_3})$ and
$\varphi(\lambda _{i_j})$, $j=1,2,3$, are linear independent over $\Z /5\Z$
 for each triple $(i_1,i_2,i_3)\in T$}\,.
\newline
From this assumption it follows that $f$ is ramified
with ramification index 5 in each $C_i$ and
each $D_{i,j,k}$.
Further, according to Lemma 1.3, $X$ is non-singular
under this assumption.
\proclaim{Lemma 2.1} Under above assumptions
\roster
\item"{($i$)}"
$C_i^2=-3$ for each $i=1,\dots,9$\,;
\item"{($ii$)}" $D_{i_1,i_2,i_3}^2=-1$ for each
$(i_1,i_2,i_3)\in T$\,;
\item"{($iii$)}"
$K_X^2=333$, where $K_X$ is the canonical class of $X$\,;
\item"{($iv$)}"
the
 geometric genera of
$C_i, 1\le i\le 9,$ and $D_{i_1,i_2,i_3},
(i_1,i_2,i_3)\in T,$ are
 $g(C_i)=4$ and $g(D_{i_1,i_2,i_3})=2$\,.
\endroster
\ep

\demo{Proof} ($i$) There are 4 triple points
of $\li$ on each $\li_i$.
Thus, $(L'_i,L'_i)=-3$.
On the other hand
$$\deg f \cdot (L'_i, L'_i)=
(f^*(L'_i),f^*(L'_i))=
(5C_i,5C_i)=25 C_i^2.$$
Therefore, $C_i^2=-3$.

Proof of ($ii$) is similar to ($i$).

($iii$) The 3-canonical class of $\widetilde{\p}^2$ is
$3K_{\widetilde{\p}^2}=-\sum L_i$.

By the pull-back formula
$$K_X=f^*
(K_{\widetilde{\p}^2})+4(\sum C_i +\sum D_{i_1,i_2,i_3}),$$
and, hence,
$$
3K_X=7\sum C_i + 12\sum D_{i_1,i_2,i_3}. \tag 2.1
$$
Thus, we have
$$
\align
9\cdot K_X^2 & =49\sum C_i^2 + 144\sum D_{i_1,i_2,i_3}^2+
168\sum
(C_i,D_{i_1,i_2,i_3})= \\
& = 49\cdot (-3)\cdot 9 +144\cdot(-1)\cdot 12
+168\cdot 4\cdot 9.
\endalign
$$
Therefore, $K_X^2=333$.

By $(2.1)$,
$$
(C_i,K_X)=9 \qquad \text{and} \qquad (D_{i_1,i_2,i_3}, K_X)=3,
\tag 2.2
$$
and  $(iv)$ follows from the adjunction formula.
\qed
\edm

\proclaim{Lemma 2.2} $X$ is a surface of general type with ample
canonical class.
\ep
\demo{Proof}
According to
the Moisheson-Nakai criterion
it is sufficient to show that
$(K_X, C)>0$ for any algebraic curve $C\subset X$. It follows from
(2.1) and (2.2) that $(K_X, C)\geq 0$ for any curve $C$. Assume that
there is an irreducible curve $C$ such that $(K_X, C)=0$.
Then the intersection of $C$ and the effective divisor
$3K_X=7\sum C_i + 12\sum D_{i_1,i_2,i_3}$ is empty. Therefore
the curve $\sigma (f(C))$ doesn't meet
any
line $L_i$, $i=1,\dots ,9$.
But it is impossible.
\qed
\edm

\proclaim{Lemma 2.3} The Euler characteristic
$e(X)$ of $X$ is
equal to $111$, and, in particular, it satisfies
the relation $K_X^2=3e(X)$.
\ep

\demo{Proof}
From $e(\widetilde{\p}^2)=15$ and $e(L_i)=e(E_{i_1,i_2,i_3})=2$
we deduce, by additivity of the Euler characteristic, that
$$\align
e(X) & =25e(\widetilde{\p}^2\setminus (\cup C_i\cup D_{i_1,i_2,i_3}))+
5\sum e(C_i\setminus \cup D_{i_1,i_2,i_3})+ \\
& + 5\sum e(D_{i_1,i_2,i_3}\setminus \cup C_i)+ 
\sum (C_i,D_{i_1,i_2,i_3}) = \\
 & =25(15-9\cdot 2-12\cdot 2+ 9\cdot 4) +5\cdot 9(2-4) +
 5\cdot 12(2-3)+9\cdot 4=111.
\endalign
$$
The relation $K_X^2=3e(X)$ now follows from Lemma 2.1 (iii). \qed
\edm
\proclaim{Corollary 2.1}
$X$ is a
strongly
rigid surface
(i.e., a surface whose
moduli space reduces to $X$ and $\bar X$ or merely to $X$,
where
$\bar X$
stands
for the complex conjugated surface).
\qed
\ep

\subheading{\S 3. Automorphisms of the coverings\ }

Let $f:X\to \widetilde{\p}^2$ be a
$(\Z /5\Z )^2$-Galois covering
considered in \S 2. Denote by $\Kl$ the group of holomorphic and
anti-holomorphic diffeomorphisms $X\to X$. Clearly, if $\Kl$ contains
at least one anti-holomorphic element, the holomorphic elements
form in $\Kl$ a subgroup $\Aut$ of index $2$. In other words, there
is a short exact sequence $1\to\Aut\to\Kl\to H\to 1$, where $H=\Z/2$
or $0$. We denote by $\kl:\Kl\to H$
the homomorphism of this sequence.
Recall that, by definition, a real structure is an anti-holomorphic
involution and note
that $H$ can be nontrivial even for varieties without real structure.

The group $\Kl$ acts most naturally on $X\times\bar X$,
$X\+ \bar X$ ($\bar X$ is the surface complex conjugated to $X$),
and the associated to them groups like $\Div$, $\Pic$,
and $H^*$, as well as on $\C(X\times\bar X)$ and $\C(X\+ \bar X)$
(where the latter is not a field, since $X\+\bar X$ is
not reducible).
There are different ways to extract from these actions
an action of $\Kl$ extending the action of $\Aut$ on
$\C(X)$, $\Div (X)$, $\Pic (X)$, and $H^*(X)$.
We choose the one which better fits to the needs of
the present investigation. In addition, it is the one
traditionally used in algebraic geometry.

To extend the action of $\Aut (X)$ on
$\C(X)$ to that of $\Kl(X)$,
we associate
with an anti-holomorphic diffeomorphism $h$ the
$\C$-anti-linear map
$h^\* :\C(X)\to\C(X)$ defined by
$h^\*(f)(x)=\overline {f(h(x))}$,
$f\in\C(X)$, $x\in X(\C)$.
The action $h^\*$ on holomorphic differential forms
is defined in a way that
$$
h^\*(df)=dh^\*(f).
$$
An anti-holomorphic diffeomorphism $h
$ defines an action on $\Div (X)$\, :
if $C\in\Div(X)$ is given by local equations
$(U_\alpha, f_\alpha)$, then $h(C)$
is given by $(h^{-1}(U_\alpha ), \overline{f_\alpha\circ h})$.
We have
$$
h^{-1}(\div f)=\div h^\*(f), f\in\C(X)
. \tag 3.1
$$
According to (3.1), $h:\Div(X)\to\Div(X)$ induces
an action $h^\* :\Pic(X)\to\Pic(X)$.
Clearly, the canonical
class $K_X\in\Pic(X)$ is invariant under
$h^\*$ for any $h\in\Kl$;
here, and further, we put $h^\*=h^*$
for $h\in\text{Aut}\, X$.
The intersection number is also preserved by any $h\in\Kl$
(it is may be worth noting that the action on
the Neron-Severi subgroup of
$H^*(X)$ associated with $h^\*:\Pic(X)\to\Pic(X)$
is not the usual $h^*$ but $-h^*$, if $h\in\Kl\setminus\Aut$).

We say that $h\in \Kl\, X$ is lifted from
$\widetilde{\p}^2$ if there exists
$\widetilde{h}\in \Kl\, \widetilde{\p}^2$
such that the following diagram is commutative
$$
\CD
X  @>{h}>> X \\
@V{f}VV            @VV{f}V     \\
\widetilde{\p ^2} @>>{\widetilde{h}}>  \widetilde{\p ^2}\, .
\endCD
$$

\proclaim{Proposition 3.1}
Every $h\in\Kl (X)$ is lifted from $\widetilde{\p}^2$.
In particular, if $X$ has a real structure then
for a proper chosen real structure of
$\widetilde{\p}^2$ the covering $f$ is defined over $\R$.
\ep

\proclaim{Lemma 3.1} Let $h\in\Kl (X)$. Then $h$ leaves fixed
the sets $\cup C_i$ and $\cup D_{i_1,i_2,i_3}$.
\ep

\demo{Proof} Assume that $h(C_{i_0})\not\subset \cup C_i$ for some $i_0$.
Then
$$
(h(C_{i_0}),\sum C_i)=a, a\geq 0.$$
It follows from the difference of genera
$g(C_{i_0})\neq g(D_{i_1,i_2,i_3})$ that
$h(C_{i_0})\neq D_{i_1,i_2,i_3}$. Therefore,
$$(h(C_{i_0}),\sum D_{i_1,i_2,i_3})=b, b\geq 0.$$
Since $h^\*(K_X)=K_X$, then by Lemma 2.1
and the adjunction formula,
$$(h(C_{i_0}),K_X)=(C_{i_0},K_X)=9.$$
Thus, in accordance with $(2.1)$ and $(2.2)$, we should have
$$7a+12b=27$$
for some non-negative integers $a$ and $b$, which is impossible.

The proof that $h(D_{i_1,i_2,i_3})\subset \bigcup_{(i_1,i_2,i_3)\in T}
D_{i_1,i_2,i_3}$ for every $(i_1,i_2,i_3)\in T$ is similar.
\qed
\edm

\demo{Proof of Proposition 3.1}
The second statement is a straightforward consequence of the first one.
To prove the latter it is sufficient to show that
$h$ acts on the fibers of $f$, i.e.,
that for almost any $p\in \widetilde{\p}^2$
one can find $q\in \widetilde{\p}^2$
such that $h(f^{-1}(p))=f^{-1}(q)$.

Let us fix a point $p_{i_0,j_0,k_0}\in \p ^2$.
Since $C_{i_0}$ and $C_{j_0}$ meet
$D_{i_0,j_0,k_0}$, then $h(C_{i_0})$ and $h(C_{j_0})$ meet
$h(D_{i_0,j_0,k_0})$. The curve $C_{i_0}$ (respectively, $C_{j_0}$)
intersects 3 other curves $D_{i_r,j_r,k_r}$, $r=1,2,3$,
(respectively, $D^{\prime}_{i_r,j_r,k_r}$, $r=1,2,3$)
distinct from $D_{i_0,j_0,k_0}$.
Thus, $h(C_{i_0})$ (respectively, $h(C_{j_0})$) intersects
each of $h(D_{i_r,j_r,k_r})$, $r=1,2,3$,
(respectively, $h(D^{\prime}_{i_r,j_r,k_r})$, $r=1,2,3$).

By Lemma 3.1,
$h(C_{i_0})=C_i$ and $h(C_{j_0})=C_j$ for some $i$ and $j$.
We have
$$\div f^*(l_{i_0}l^{-1}_{j_0})=
5(C_{i_0}+\sum_{r=1}^{3}D_{i_r,j_r,k_r})-5(C_{j_0}+\sum_{r=1}^{3}
D^{\prime}_{i_r,j_r,k_r})$$
and
$$
\div f^*(l_{i}l^{-1}_{j})=
5(h(C_{i_0})+\sum_{r=1}^{3}h(D_{i_r,j_r,k_r}))-5(h(C_{j_0})+\sum_{r=1}^{3}
h(D^{\prime}_{i_r,j_r,k_r})).
$$
Therefore, there is a constant $k_{i_0,j_0}$ such that
$$ h^\*(f^*(l_{i}l^{-1}_{j}))=k_{i_0,j_0}f^*(
l_{i_0}l^{-1}_{j_0}).
\tag 3.2$$

Let us
choose another point $p^{\prime}_{i_0,j_0,k_0}\in \p ^2$,
$p_{i^{\prime}_0,j^{\prime}_0,k^{\prime}_0}\in \li _{i'_0}\cap \li _{j'_0}$
and consider the curves $C_{i^{\prime}_0}$,
$C_{j'_0}$, and their images $h(C_{i^{\prime}_0})=C_{i^{\prime}}$,
$h(C_{j^{\prime}_0})=C_{j^{\prime}}$.
Arguing as above, we conclude
that there exists a constant $k_{i'_0, j'_0}$ such that
$$h^\*(f^*(l_{i^{\prime}}l^{-1}_{j^{\prime}}))=
k_{i'_0,j'_0}f^*(l_{i^{\prime}_0}l^{-1}_{j^{\prime}_0}).
\tag 3.3  $$
Since
every $p\in \widetilde{\p}^2\setminus \cup D_{i_1,i_2,i_3}$ can be given
as the intersection of fibers of two linear rational functions
$l_{i_0}l^{-1}_{j_0}$ and $l_{i^{\prime}_0}l^{-1}_{j^{\prime}_0}$,
it follows from $(3.2)$ and $(3.3)$ that for any
$p\in \widetilde{\p}^2\setminus \cup D_{i_1,i_2,i_3}$ we have
$h(f^{-1}(p))=f^{-1}(q)$ for some $q\in \widetilde{\p}^2$.
\qed
\edm

\subheading{\S 4.\ Three examples}

{\bf Example I.} {\it 
A non real rigid surface.}

Let $\li=\li _1\cup \dots \cup \li _9$ be  an arrangement of nine
lines in $\p ^2$ dual to the nine inflection points
of a smooth cubic $C$ in the dual plane (see \S 2),
and let $f: X_1\to \widetilde{\p}^2$ be
the Galois covering associated with the epimorphism
$\ov{\varphi}_1 : H_1(\p ^2 \setminus \li ,\Z )\to (\Z /
5\Z)^2$ given by
$$
\matrix \ov{\varphi}_1(\lambda _1)=(1,1), & \ov{\varphi}_1(\lambda _2)=
(1,0), & \ov{\varphi}_1(\lambda _3)=(1,1), \\
\ov{\varphi}_1(\lambda _4)=(3,3), & \ov{\varphi}_1(\lambda _5)=
(3,0), & \ov{\varphi}_1(\lambda _6)=(0,1), \\
\ov{\varphi}_1(\lambda _7)=(0,1), & \ov{\varphi}_1(\lambda _8)=
(0,2), & \ov{\varphi}_1(\lambda _9)=(1,1),
\endmatrix
$$
see \S 1 (note that $\sum\varphi_1(\lambda_i)=0\mod 5$).

\proclaim{Proposition 4.1} The surface $X_1$ is smooth and
strongly rigid. The group $\Kl(X_1)$ coincides with the
covering transformations group $G=\Z/5\times\Z/5$.
In particular, there does not exist
neither a real structure nor even an
anti-holomorphic diffeomorphism on $X_1$.
\ep

\demo{Proof}
The surface $X_1$
is smooth due to Lemma 1.3.
According to Lemmas 2.1, 2.3 we have
$K_{X_1}^2=333$ and $e(X_1)=111$, and the rigidity
statement follows from Corollary 2.1.

Consider, now,
any $c\in\Kl(X_1)$.
By Proposition 3.1, $c$ is lifted from $\widetilde{\p}^2$, i.e.,
there is $\widetilde c\in \Kl(\widetilde{\p}^2)$
such that $f\circ c=\widetilde{c}\circ f.$

As in \S 1, consider affine coordinates $x_1,x_2$ in
$\C ^2=\p ^2\setminus \li _{\infty}$ and the linear equations
$l_i(x_1,x_2)=0$ of $\li _i \cap \C ^2$.
The function field $\C (X_1)$ of $X_1$ is identified
with the sub-field
$$K_{\ov{\varphi_1}}=\C (x_1,x_2,w_1,w_{2})$$ of $K_{u(5)}$,
where $w_1=l_1l_2l_3l_4^3l_5^3l_9$ and
$w_2=l_1l_3l_4^3l_6l_7l_8^2l_9$,
so that
$$K_{\ov{\varphi_1}}=\bigoplus_{a\in \pr A_1}\C (x_1,x_2)
z^{a}
\tag 4.1$$
is a subspace of the vector space
$$K_{u(m)}=\bigoplus_{a\in \pr A}\C (x_1,x_2)z^a
$$
over $\C (x_1,x_2)$, where
$$A=\{
\alpha=(\alpha_1,\dots,\alpha_9)\in \Z^9 \, \mid \,
0\le \alpha_i\le 4 \, \text{and} \, \sum
\alpha_i=0\mod 5 \} ,$$
$\text{pr} : A\mapsto \overline{A}\simeq {(\Z/5\Z)}^8$ is the projection
given by
$\text{pr}(\alpha )=
(\alpha _1,\dots ,\alpha _8)$ for
$\alpha=(\alpha _1,\dots ,\alpha _9)$, and
$A_1\subset A$
consists from $0=(0,0,0,0,0,0,0,0,0)$ and
$$
\matrix
(1,1,1,3,3,0,0,0,1), & (2,2,2,1,1,0,0,0,2), & (3,3,3,4,4,0,0,0,3), &
(4,4,4,2,2,0,0,0,4),\\
(1,0,1,3,0,1,1,2,1),  & (2,0,2,1,0,2,2,4,2), & (3,0,3,4,0,3,3,1,3), &
(4,0,4,2,0,4,4,3,4), \\
(2,1,2,1,3,1,1,2,2), & (4,2,4,2,1,2,2,4,4), & (1,3,1,3,4,3,3,1,1), &
(3,4,3,4,2,4,4,3,3), \\
(3,1,3,4,3,2,2,4,3), &  (1,2,1,3,1,4,4,3,1), & (4,3,4,2,4,1,1,2,4),  &
(2,4,2,1,2,3,3,1,2), \\
(4,1,4,2,3,3,3,1,4), & (3,2,3,4,1,1,1,2,3), & (2,3,2,1,4,4,4,3,2), &
(1,4,1,3,2,2,2,4,1), \\
(0,1,0,0,3,4,4,3,0), & (0,2,0,0,1,3,3,1,0), &  (0,3,0,0,4,2,2,4,0), &
(0,4,0,0,2,1,1,2,0).
\endmatrix
$$

The diffeomorphism
$c$ induces an action
$c^{\*}$ on $\C (X_1)$
such that the restriction of $c^\* $ to the subfield $\C (\widetilde{\p}^2)=
\C ({\p}^2)$ coincides with $\widetilde c^\*$ (see Section 3).
By Lemma 3.1, the sets $\cup C_i$ and $\cup D_{i_1,i_2,i_3}$
are invariant under the action of $c$. Therefore the set $\cup \li_i$ is
invariant under the action of $\widetilde c$.
Thus, $c^\*$ acts on the set of the one-dimensional
subspaces
$\C (x_1,x_2)z^a$, $\ov{a}\in \pr A_1$, of
$K_{\ov{\varphi_1}}$,
and, thus, induces an action on $A_1$.
We denote the latter action by $c^\*$ also.
For $a\in A_1$ denote by $r_i(\ov{a}), i\in \Z/5\Z,$ the number
of coordinates of $\ov{a}$ equal $i$.

\proclaim{Lemma 4.1} The functions $r_i$ are
invariant under the action of $c^\*$, i.e.,
$r_i(\alpha )=r_i(\beta )$ for $\beta =c^\* (\alpha )$.
\ep
\demo{Proof}
For each $j$, $1\le j\le 9$
the coordinate $\alpha _j$ of
$\alpha=(\alpha _1,\dots ,\alpha _9) \in A_1$
is congruent modulo 5 to the order of zero along $C_j$
of any of the functions
in $C(x_1,x_2)z^a$, $a={\pr \alpha}$.
It remains to note that due to Lemma 3.1 $c$ interchanges the curves
$C_j$.
\qed
\edm
By Lemma 4.1, the action of $c^\*$ on $A_1$ is determined by a
permutation $\pi$ of $1,\dots, 9$.

Consider
$\alpha
=(1,1,1,3,3,0,0,0,1)$ and
$\beta =(1,0,1,3,0,1,1,2,1)$.
It is easy to see that $\alpha $ is
the unique element in $A_1$ with
$r_0=3,\, r_1=4,\, r_2=0,\, r_3=2, r_4=0$.
Respectively,
$\beta $ is the unique element in $A_1$ with
$r_0=2,\, r_1=5,\, r_2=1,\, r_3=1, r_4=0$.
Thus,  by Lemma 4.1,
$c^\* (\alpha )=\alpha $ and $c^\* (\beta )=\beta $.
Since $r_2(\beta )=1$ and $r_3(\beta )=1$,
we have $\widetilde{c}(\li_4)=\li_4$ and $\widetilde{c}(\li_8)=\li_8$.
Further, $r_3(\alpha )=2$ implies
$\widetilde{c}(\li_5)=\li_5$
and $r_0(\beta )=2$ implies
$\widetilde{c}(\li _2)=\li _2$.

The above invariance properties of
$\li _2$, $\li _4$, $\li _5$, $\li _8$
mean that these lines
are
invariant
under the action of $\widetilde c$.
Hence, the points
$p_{2,4,9}=\li_2\cap\li_4$, $p_{2,5,8}=\li_5\cap\li_8$, $p_{4,5,6}
=\li_4\cap\li_5$ and $p_{3,4,8}=\li_4\cap\li_8$
are fixed points of $\widetilde c$.

Since $r_0(\alpha)=3$,
there remain two possibilities: either
$\widetilde{c}(\li _6)=\li _7$ and $\widetilde{c}(\li _7)=\li _6$, or
$\widetilde{c}(\li _6)=\li _6$ and $\widetilde{c}(\li _7)=\li _7$.

If $\widetilde{c}(\li _6)=\li _7$ and $\widetilde{c}(\li _7)=\li _6$, then
their intersection point $p_{2,6,7}$
is a fixed point.
This is impossible. Indeed, in this case $\li _6$
passes through two different fixed
points $p_{2,6,7}$ and
$p_{4,5,6}$, so should
satisfy $\widetilde{c}(L_6)=L_6$.

If $\li _6$ and $\li _7$ are invariant
lines, then all
lines $\li _i$ with $1\le i\le 9$, should be
invariant.
In fact, since $\li _5$ and $\li _7$ are invariant
lines,
their intersection point $p_{3,5,7}$ is a fixed point.
Then, $\li _3$ is an invariant
line, since $\li _3$ passes through
two fixed
points $p_{3,4,8}$ and $p_{3,5,7}$. Therefore,
the intersection points $p_{3,6,9}$ of
$\li _3$ and $\li _6$, $p_{1,2,3}$ of
$\li _2$ and $\li _3$, $p_{1,4,7}$ of
$\li _4$ and $\li _7$, and $p_{7,8,9}$ of
$\li _7$ and $\li _8$ are also fixed points.
It implies, that $\li _1$ and $\li _9$,
which go, respectively, through
$p_{1,2,3}$, $p_{1,4,7}$ and $p_{3,6,9}$, $p_{7,8,9}$
are invariant under the action of $\widetilde{c}$.
0
As we have proved, the nine
inflection points
of $C$, which is a smooth cubic,
are
fixed
under the action induced on $\p^2$ by
$\widetilde{c}$.
Therefore, if $\widetilde{c}\in\Aut(\widetilde\p^2)$,
then $\widetilde{c}=\Id$ and, hence, $c$ is a covering transformation.
If  $\widetilde{c}\notin\Aut(\widetilde\p^2)$,
then $\widetilde{c}^2\in\Aut{\widetilde\p^2}$
is the indentity, and so
$\widetilde{c}$ induces a real structure on $\p ^2$
such that all inflection points of a smooth cubic $C$ are real with respect
to this structure, but it is impossible. \qed
\edm

\proclaim{Corollary 4.1}
The moduli space of complex structures on
the underlying smooth 4-manifold consists of
two distinct points $X_1$ and $\bar X_1$.
In particular, $X_1$ and $\bar X_1$ give a counterexample to
"Dif=Def" problem
\footnote{First counter examples to "Dif=Def" problem were given by
Manetti in \cite{Ma}.}.
\ep

{\bf Remark 4.1.}
One can deduce from Proposition 4.1 and
the Mostow strong
rigidity (using the Smith inequality for transformations
of prime order and group cohomology arguments) that
$X_1$ has no
any
nonindentical
diffeomorphism of order $\ne 5$.

{\bf Remark 4.2.}
The irregularity of $X_1$ is equal to zero.
It follows, for example, from \cite{Is}.

{\bf Example II.} {\it A non maximal rigid surface.}

Now, suppose that the cubic $C$ is given
by $x_1^3+x_2^3+x_3^3=0$ and that the lines
$\li_1,\dots,\li_9$ are numbered as in Section 2.
In particular, under this choice,
$\li_1, \li_5$, and $\li_9$
are real.
Let $f: X_2\to \widetilde{\p}^2$ be the
Galois covering associated with the epimorphism
$\ov{\varphi}_2 : H_1(\p ^2 \setminus \li ,\Z )\to (\Z /m\Z)^2$ given by
$\ov{\varphi}_2(\lambda _i)=(a_{i,1},a_{i,2})$, where
$$
\matrix \ov{\varphi}_2(\lambda _1)=(0,1), & \ov{\varphi}_2(\lambda _2)=
(1,0), & \ov{\varphi}_2(\lambda _3)=(1,0), \\
\ov{\varphi}_2(\lambda _4)=(0,1), & \ov{\varphi}_2(\lambda _5)=
(1,0), & \ov{\varphi}_2(\lambda _6)=(0,1), \\
\ov{\varphi}_2(\lambda _7)=(1,2), & \ov{\varphi}_2(\lambda _8)=
(1,2), & \ov{\varphi}_2(\lambda _9)=(0,3).
\endmatrix
$$

\proclaim{Proposition 4.2} The surface $X_2$ is smooth and
strongly rigid.
It can  be equipped with a real structure. Such a structure
is unique,
up to conjugation by covering transformations,
and not maximal
{\rom
(where the latter means that $\sum \dim H_i(X_2(\R );$ $\Z/2\Z)<
\sum\dim H_i(X_2(\C);\Z/2\Z)).$ The group $\Kl(X_2)$ is
a semi-direct product of the group $\mu _2\simeq \Z/2$ of order 2
and  the covering transformations group
$G\simeq \Z/5\times \Z/5$. This
$\Z/2$-extension is defined by relations
$s\gamma s^{-1}=
\gamma^{-1}$, $\gamma\in G$, $s\in \mu_2$, $s\neq 1$.
\ep

\demo{Proof}
As in the proof of Proposition 4.1,
$X_2$ is smooth due to Lemma 1.3.
According to Lemmas 2.1, 2.3 we have
$K_{X_2}^2=333$ and $e(X_2)=111$, and the rigidity
statement follows from Corollary 2.1.

As above, we identify the function field $\C (X_2)$
of $X_2$ with subfield
$K_{\ov{\varphi_2}}=\C (x_1,x_2,w_1,w_{2})$ of $K_{u(5)}$, where
$w_1=l_2l_3l_5l_7l_8$ and $w_2=l_1l_4l_6l_7^2l_8^2l_9^3$.
Then
$$K_{\ov{\varphi_2}}=\bigoplus_{a\in \pr A_2}
\C (x_1,x_2)z^a
$$
is a subspace of the vector space
$$K_{u(m)}=\bigoplus_{a\in \pr A}\C (x_1,x_2)z^a
$$
over $\C (x_1,x_2)$, where
$A$ and $\pr$ are the same as in the previous example
and where
$A_2$
consists from $(0,0,0,0,0,0,0,0,0)$
and
$$
\matrix
(0,1,1,0,1,0,1,1,0), & (0,2,2,0,2,0,2,2,0), & (0,3,3,0,3,0,3,3,0), &
(0,4,4,0,4,0,4,4,0),\\
(1,0,0,1,0,1,2,2,3),  & (2,0,0,2,0,2,4,4,1), & (3,0,0,3,0,3,1,1,4), &
(4,0,0,4,0,4,3,3,2), \\
(1,1,1,1,1,1,3,3,3), & (2,2,2,2,2,2,1,1,1), & (3,3,3,3,3,3,4,4,4), &
(4,4,4,4,4,4,2,2,2), \\
(1,2,2,1,2,1,4,4,3), &  (2,4,4,2,4,2,3,3,1), & (3,1,1,3,1,3,2,2,4),  &
(4,3,3,4,3,4,1,1,2), \\
(1,3,3,1,3,1,0,0,3), & (2,1,1,2,1,2,0,0,1), & (3,4,4,3,4,3,0,0,4), &
(4,2,2,4,2,4,0,0,2), \\
(1,4,4,1,4,1,1,1,3), & (2,3,3,2,3,2,2,2,1), &  (3,2,2,3,2,3,3,3,4), &
(4,1,1,4,1,4,4,4,2).
\endmatrix
$$

Pose $\alpha =(0,1,1,0,1,0,1,1,0)$ and $\beta =(1,0,0,1,0,1,2,2,3)$
and consider any $c\in \text{Kl}(X_2)$, $c\ne\Id$.
The considerations as in the proof of Proposition 4.1
show that
the line $\li_9$ and
each
the unions $\li_7\cup\li_8$, $\li_1\cup\li_4\cup\li_6$,
and $\li_2\cup\li_3\cup\li_5$ are invariant under the action
of $\widetilde{c}$.

It is impossible that
$\widetilde{c}(\li _7)=\li _7$ and
$\widetilde{c}(\li _8)=\li _8$.
In fact,
otherwise, $\widetilde{c} (p_{1,4,7})=p_{1,4,7}$, since
the arrangement $\li_1\cup\li_4\cup\li_6\cup\li_7\cup\li_8$
has only two $3$-fold points $p_{1,4,7}$ and $p_{1,6,8}$.
It would imply $\widetilde{c}(\li_6)=\li_6$, which together
with $\widetilde{c}(\li_9)=\li_9$ implies that $\li_1$, and hence
$\li_4$ and subsequantly all the lines, are invariant under
the action of $\widetilde{c}$, which contredicts to $\widetilde{c}\ne\Id$.

So, the only possibility is $\widetilde{c}(\li_7)=\li_8$
and $\widetilde{c}(\li_8)=\li_7$.
Since the pair of the $3$-fold
points $\{p_{2,5,8},p_{3,5,7}\}$
of $\li_2\cup\li_3\cup\li_5\cup\li_7\cup\li_8$
is invariant under $\widetilde{c}$, the line $\li_5$
is invariant while $\li_2$ and $\li_3$
are permuted.
The same arguments show that
$\widetilde{c}(\li_1)=\li_1$,
$\widetilde{c}(\li _4)=\li _6$, and $\widetilde{c}(\li_6)=\li_4$.

Such an action of $\widetilde{c}$ on $L=\cup \li_i$
is
the one induced by the standard complex conjugation
on $\widetilde{\p}^2$ (see Section 2 or use the unicity) and
thus coincides with it.
It lifts to a real structure $s$ on $X_2$; in fact,
$X_2$
can be seen as the minimal desingularization
of the projective closure of the
real surface given by equations
$$
\align
 w_1^5  = & (x_1^2+x_1+1)(x_2-1)(x_1^2+x_1x_2+x_2^2), \\
 w_2^5  = & (x_1-1)(x_2^2+x_2+1)(x_1^2+x_1x_2+x_2^2)^2(x_1-x_2)^3.
\endalign
$$
This real surface
is not maximal, since its real part is homeomorphic
to $\RP $ with four blown up points (it is easy to check that there are
only four real points among the
blown up
points
$p_{i,j,k}$, $\{i,j,k\}\in T$).
Since each $c\in\Kl(X_2)$
is determined, up to composition with covering
transformations, by
$\widetilde{c}$,
the group
$\Kl(X_2)$ is generated by $s$
and the covering transformations.
The commutation relations
$s\gamma=\gamma^{-1} s$
follow from the above equations.
These relations  imply that each
$s\gamma$ is a real structure and that these
real structures are all equivalent.
\qed\edm
\bk

{\bf Remark 4.3.}
The surfaces in the both examples
considered have the same
$K^2$ and $e$. Thus, they belong to the same Hilbert scheme
and provide an example of a Hilbert scheme
whose connected components have
different properties with respect to the existence of real structures
on the surfaces
repersenting these components. Note also that
contrary to the first example in the second one
the moduli space reduces  to one point, which is real
(and, moreover, corresponds to a surface with a real structure).

{\bf Example III.} {\it 
A rigid surface with two non-equivalent real structures.}

Here,
we call
two structures equivalent
if they can be transformed one into another by
an automorphism of the surface.

Let
$L=L_1\cup \cdots \cup L_6$ be
a complete quadrilateral. Note that
two complete quadrilaterals are projectively equivalent.
In this arrangement
$t_2=3$, $t_3=4$, and $t_r=0$ for $r\geq 4$.
After suitable numbering, we can assume that the set of $2$-fold points
consists of $\{L_1\cap L_4,\, L_2\cap L_5,\, L_3\cap L_6 \}$ and
the set of $3$-fold points
of $\{L_1\cap L_2\cap L_6,\, L_2\cap L_3\cap L_4,\,
L_1\cap L_3\cap L_6,\, L_4\cap L_5\cap L_6 \}.$

Let $f: X_3\to \widetilde{\p}^2$ be
the Galois covering associated with the epimorphism
$\ov{\varphi}_3 : H_1(\p ^2 \setminus \li ,\Z )\to (\Z /
5\Z)^2$ given by
$\ov{\varphi}_3(\lambda _i)=(a_{i,1},a_{i,2})$, where
$$
\matrix \ov{\varphi}_3(\lambda _1)=(1,0), & \ov{\varphi}_3(\lambda _2)=
(1,0), & \ov{\varphi}_3(\lambda _3)=(1,2), \\
\ov{\varphi}_3(\lambda _4)=(0,1), & \ov{\varphi}_3(\lambda _5)=
(0,1), & \ov{\varphi}_3(\lambda _6)=(2,1),
\endmatrix
$$
and
$\widetilde{\p}^2$ is the blow up of $\p ^2$ at the $3$-fold points
of $L$.
As above, denote by $\sigma :\widetilde {\p ^2} \to \p ^2$ the blow up
with the centers at
the $3$-fold points,
by $E_{i,j,k}$ the exceptional divisor over the $3$-fold point $p_{i,j,k}$,
and by $L'_i$ the strict transform
of $L_i$.
Pose $C_i=f^{-1}(L'_i)$ and $D_{i,j,k}=f^{-1}(E_{i,j,k})$.

As in \S 1, consider affine coordinates $x_1,x_2$ in
$\C ^2=\p ^2\setminus L _{\infty}$ and the linear equations
$l_i(x_1,x_2)=0$ of $L_i \cap \C ^2$. Then, by Lemma 1.3, $X_3$ is
isomorphic to the minimal desingularization
of the projective closure of the surface given by equations
$$
\align
 w_1^5  = & l_1l_2l_3l_6^2, \tag 4.2 \\
 w_2^5  = & l_3^2l_4l_5l_6. \tag 4.3
\endalign
$$

The computations as in the proof of Lemmas 2.1 and 2.2 show that $X_3$
is a surface of general type with
$K_{X_3}^2=45$ and $e(X_3)=15$. Therefore $X_3$ is a
strongly rigid surface.
\proclaim{Lemma 4.3}
 Let $h\in\Kl (X_3)$. Then $h$ leaves fixed
the set $(\cup C_i)\cup (\cup D_{i_1,i_2,i_3})$.
\ep
\noindent
{\it Proof} is similar to that of Lemma 3.1.
\proclaim{Proposition 4.3}
Every $h\in\Kl (X_3)$ is lifted from $\widetilde{\p}^2$.
In particular, if $X_3$ has a real structure then
for a proper chosen real structure of
$\widetilde{\p}^2$ the covering $f$ is defined over $\R$.
\ep
\noindent
{\it Proof} is similar to that of Proposition 3.1.
\proclaim{Proposition 4.4} The surface $X_3$ can
be equipped with 2 non-equivalent real structures.
\ep
\demo{Proof}
Consider two real structures of $\p ^2$. For the first one, all
the
lines
of $L$ are real, and for the second one, the lines $L_3$, $L_6$
are real and the lines $L_1$, $L_2$, respectively $L_4$ and $L_5$, are
complex conjugated. Then these two real structures induce two real
structures on $X_3$, since in
the
both cases
the polynomials in
$(4.2)$ and $(4.3)$
are
defined over $\R$.

These two real structures of $X_3$ are non-equivalent.
Indeed,
by Lemma 4.3,
each automorphism
of $X_3$
leaves fixed
the set $(\cup C_i)\cup (\cup D_{i_1,i_2,i_3})$
while, on one hand,
all the curves $C_i$
and
$D_{i_1,i_2,i_3}$ are real with respect to
the first real structure, but,
on the other hand,
only $C_3$ and $C_6$ (among
$C_1\dots , C_6$) are
real curves with respect to
the second real structure.
\edm

\subheading{\S 5.\ Non reality of fake projective plane and remarks}

{\bf A.}
We call a surface of general type with
$p_g=q=0$ and $K^2=9$ {\it a fake projective plane}.
The existence of fake projective planes was proved by D. Mumford
\cite{Mu}.

\proclaim{Theorem 5.1} A fake
projective plane has no anti-holomorphic
diffeomorphisms.
\ep
\demo{Proof} Let $X$ be a fake projective plane. Then
(see \cite{Mi}, \cite{Y}),
the universal covering space of $X$ is a ball.

First, let us show
that there is no an anti-holomorphic involution
on $X$.

So, assume that $X$ can be equipped with a real structure
and denote by $X_\R$ the real point set of $X$.
According to the Lefschetz trace formula
applied to the involution defining the structure,
$e(X_\R)=1$. Thus, $X_\R$ is in nonempty and contains at least
one component diffeomorphic either to
sphere or
real projective plane.
To lift the real structure to the universal covering pick a point $p$
on such a component and identify the points of the universal covering with
homotopy classes of the paths starting at $p$.
The real part of the covering covers (without ramification) the chosen
real component of $X$. On the other hand, since the universal covering
space is a ball, its real part has no compact components.

From the above it follows now
that if there
exists
an anti-holomrphic
diffeomorphism $h$, then its
order
can not be $2n$, where $n$ is odd. In fact,
if $n$ is odd, then $h^n$ is the anti-holomorphic involution.
Thus, Theorem 5.1
follows from the following Lemma.

\proclaim{Lemma 5.1}
The group $\text{Aut}\, X$
have no elements of even order.
\endproclaim

\demo{Proof of Lemma} Assume that there is $h\in \text{Aut}\, X$ of
order 2.
One dimensional components $C$ of the fixed point set of $h$
are nonsingular. By
Enoki-Hirzebruch \cite{BHH}
relative proportionality,
$$ e(C)=2C^2. $$
Thus $C=\emptyset $, since
otherwise $C^2>0$ and $e(C)<0$
(the latter inequalities can be deduced, for example,
from $C=rK$ with positive $r\in\Q$).

Since $\dim H^i(X,\C)=1$ for $i=0,2,4$ and $0$ for $i=1,3$,
the topological Lefschetz trace formula
shows that
the number of fixed points of $h$
should be equal to 3 for any nontrivial
holomorphic automorphism
without one dimensional components in the fixed point
set.
Next, applying
the holomorphic Lefschetz formula
to such a $h$ (of order 2), we
get
$$\sum _{i=1}^3\frac{1}{\det (Id-D_i)}=1,$$
where $D_i, i=1,2,3,$
are the Jacobi
matrices
of $h$ at its
fixed points.
On the other hand,
$\det (Id-D_i)=4$
at each fixed point, and thus
the Lefschetz formula turns into
$\frac{3}{4}=1$, i.e., we get a contradiction,
which proves the Lemma and finishes the proof of Theorem 5.1.
\qed
\enddemo

\proclaim{Corollary 5.1}
For any fake projective plane $X$
the moduli space of complex structures on
the underlying smooth 4-manifold consists of
two distinct points $X$ and $\bar X$.
\ep

\edm

{\bf B.} Arguments used in the proof of Theorem 5.1
to exclude anti-holomorphic involutions
can be replaced by the following general result.

\proclaim{Theorem 5.2}
If
$X$
is a compact complex K\"ahler surface
of negative sectional curvature,
then
for any real structure
the real part of $X$ has no real component diffeomorphic to
sphere, real projective plane, torus or Klein bottle.
\ep
\demo{Proof} Let $p:B\to X$ be the universal covering.
According to Cartan-Hadamard theorem, $B$ is diffeomorphic
to $\R^4$.
Each connected component $M$ of the pull back $p^{-1}(F)$
of a real component $F$ of $X$ is a real component
of some real structure on $B$. Hence, by Smith
theorem,
$M$ has the homology of a point and, hence,
it is diffeomorphic to $\R^2$.
This excludes
sphere and real projective plane as $F$ and gives the injectivity
of $\pi_1(F)\to\pi_1(X)$. It remains to note that $\pi_1(X)$,
as
the fundamental group of a compact
manifold of negative curvature,
contains no subgroup
isomorphic to $\Z\oplus\Z$, see \cite{P}.
\qed
\enddemo

It may be interesting to compare this observation with
Koll\`ar conjecture (and Viterbo theorem, see \cite{Kh})
according to which an algebraic variety
of dimension $\ge 3$ is of general
type as soon as one of its real components,
with respect to some real structure, is hyperbolic.

{\bf C.}
Miyaoka-Yau surfaces
can
provide interesting examples related to the "Ragsdale
bound",
i.e.,
examples of real surfaces $X$ with
$\beta_1^\R=\dim H_1(X_\R ,\Z/2\Z)$ close or above $h^{1,1}(X)$.
(First examples with $\beta_1^\R>h^{1,1}(X)$ were found in
early 80th by I.~Itenberg \cite{It}.)
Recall that a real surface $X$ is called
{\it maximal} (or $M$-surface), if
the Smith bound (see, f.e., the survey
\cite{DK})
$$\sum \beta_i^\R\le\sum \beta_i^\C=
2+4(h^{1,0}+\nu)+2h^{2,0}+h^{1,1},
\tag 5.1 $$
where $\nu$ is the rank of the $2$-torsion in $H_1(X;\Z)$
and $\beta_i^\C=\dim H_i(X;\Z/2\Z)$,
turns into equality.

By the Lefschetz formula, for any real surface
$$\beta_0^\R-\beta_1^\R+\beta_2^\R=1+\tr P^{1,1},
\tag 5.2$$
where $P^{1,1}$ is the primitive part of $H^{1,1}$ (which is, in fact,
of codimension $1$ in $H^{1,1}$).
Hence, for an $M$-surface
$$\beta_1^\R=1+2(h^{1,0}+\nu)+h^{2,0}+p^{1,1}_-,
\tag 5.3 $$
where $p^{1,1}_-$ stands for the dimension of
the anti-invariant part of the action
of the real structure in $P^{1,1}$.
On the other hand, for Miyaoka-Yau surfaces
$$
3(2+2h^{2,0}-h^{1,1})=2-4h^{1,0}+2h^{2,0}+h^{1,1}$$
and thus
$$
h^{1,1}=h^{2,0}+h^{1,0}+1.
\tag 5.4 $$
Finally, {\it for any maximal real Miyaoka-Yau surface}
$$\beta_1^\R=h^{1,1}+p^{1,1}_-+h^{1,0}+2\nu. \tag 5.5 $$
It implies that either for all maximal real Miyaoka-Yau
surfaces with $h^{1,0}=0$ it holds
$p^{1,1}_-=\nu=0$ (which would be strange)
or there are (maximal) real Miyaoka-Yau surfaces with $h^{1,0}=0$
and
$\beta_1^\R>h^{1,1}$ (which is more plausible).

The next propositions
show that if maximal real Miyaoka-Yau surfaces
of general type exist their topology should
be very restricted.
Note also that
Theorem 5.2 provides a below bound on
$\vert e(F)-1\vert$ for the real components $F$ of $X$, while the more
traditional results give the upper bounds of
$\vert e(X_\R)-1\vert$
where $X_\R$ is the whole real point set
(see, f.e., the survey \cite{DK}).

\proclaim{Proposition 5.1}
There
is no any maximal real
Miyaoka-Yau surface
with $h^{2,0}\leq 3$.
\ep
\demo{Proof} Let $X$ be a maximal real Miyaoka-Yau surface.
Denote by $k$ the number of connected components of $X_\R$.
By (5.2),
$$2k -\beta_1^\R = 1 +p_+^{1,1}-p_-^{1,1}.
\tag 5.6 $$
Substituting $\beta_1^\R$ from $(5,5)$ into $(5.6)$, and via $(5.4)$,
we have
$$ 2k=h^{1,1}+h^{1,0}+2\nu +
p_+^{1,1}+1
\tag 5.7 $$
By Theorem 5.2,
$\dim H_1(S;\Z/2\Z)\geq 3$
for any connected component
$S$ of $X_\R$
(the theorem is applied,
since the unviversal covering of $X$
is a ball, see \cite{Mi},\cite{Y}). Therefore $\beta_1^\R \geq 3k$ and, by $(5.2)$ and $(5.7)$
$$2p_-^{1,1}\geq h^{1,1}+h^{1,0}+2\nu +3p_+^{1,1}+3$$
or, equivalently,
$$h^{1,1}\geq h^{1,0}+2\nu +5p_+^{1,1}+5
\tag 5.8$$
and
$$h^{2,0}\geq 2\nu +5p_+^{1,1}+4. $$
Therefore $h^{2,0}\geq 4$.
\qed
\edm

\proclaim{Proposition 5.2}
Let $X$ be a maximal real
Miyaoka-Yau surface.
Then there are at least 3 connected components of $X_\R$
diffeomorphic
to a sphere with 3 points blown up.
\ep
\demo{Proof}
Denote by $k_3$ the number of connected components of $X_\R$
diffeomorphic to a sphere with 3 points blown up.
Then, by Theorem 5.2,
$\dim H_1(S;\Z/2\Z)\ge 4$
for all other connected components
$S$ of $X_\R$. Therefore $\beta_1^\R \geq 4k- k_3$. It follows from
$(5.5)$ and $(5.7)$ that in this case we should have the following
inequality
$$h^{1,1}+h^{1,0}+2\nu +p_-^{1,1}\geq 2h^{1,1}+2h^{1,0}+4\nu +
2p_+^{1,1}+
2 -k_3$$
which
contradicts to $h^{1,1}>p_-^{1,1}$ if $k_3<3$.
\qed
\edm

{\bf D.} The Miyaoka-Yau surfaces are
quasi-simple in the following sense: two real structures
of
such a surface
are conjugated by
an automorphism,
if and only if they are conjugated by
a diffeomorphism.\footnote{In general, the real quasi-simplicity
of a deformation class of varieties should mean that
two real structures are real deformation equivalent, if and only
if they are conjugated by a diffeomorphism.}
This follows from Mostow strong rigidity
and the fact that the only isometry of a compact hyperbolic
riemannian manifold acting identically on
the fundamental group is the identity map. (Note, that
two real structures are conjugated by an element of $\Aut$
as soon as they are conjugated by an element of $\Kl$.)

{\bf Added in proof.}
Using the nonreal surface constructed in Section 4 or
fake projective planes (see Section 5),
one can obtain examples
of varieties~$X$ of any dimension $\ge3$ having the
same property, i.e.,
examples such that $X$ and $\bar X$ belong to
distinct connceted components of the moduli space.
It is sufficient to consider products of these surfaces
with tori.
The statement on the components of the moduli space will then follow from
the well known properties of the Albanese map and Siu's rigidity theorem.

F.~Catanese informed us that he also constructed
examples of surfaces
where the complex conjugation interchanges the components
of their moduli space.
His surfaces are covered by the bi-disc $D\times D\subset
\C^2$.

\end